\documentclass[11pt,oneside]{article}

\usepackage{amsthm, amsmath, amssymb, amsfonts, epsfig, mathdots}
\usepackage{mathtools}
\usepackage{esint}
\usepackage{epstopdf}

\usepackage{float,graphicx}
\usepackage[font=sl,labelfont=bf]{caption}
\usepackage{subcaption}

\usepackage{enumerate}

\usepackage[colorlinks=true, pdfstartview=FitV, linkcolor=blue,
            citecolor=blue, urlcolor=blue]{hyperref}
\usepackage[usenames]{color}
\definecolor{Red}{rgb}{0.7,0,0.1}
\definecolor{Green}{rgb}{0,0.7,0}

\usepackage[numbers]{natbib}

\usepackage{url}
\makeatletter\def\url@leostyle{%
 \@ifundefined{selectfont}{\def\UrlFont{\sf}}{\def\UrlFont{\scriptsize\ttfamily}}} \makeatother

\usepackage{accents, comment}

\usepackage{bbm}

\usepackage{mathrsfs}

\newtheorem{theorem}{Theorem}

\theoremstyle{definition}

\theoremstyle{remark}
\newtheorem{remark}[theorem]{Remark}

\numberwithin{equation}{section}
\numberwithin{theorem}{section}

\newcommand{\norm}[1]{ \| #1 \| }       

\title{On higher-order derivative ratios in turbulent flows}

\vspace{.7in}
\author{
Zoran Gruji\'c \ \ \ \ \ \ Muhammad Mohebujjaman\\
\small UAB  \\
}

\begin{document}
\maketitle

\begin{abstract}
A computational study of higher-order derivative ratios on a time interval leading to the enstrophy peak is presented in the case of the
3D Taylor-Green vortex, a benchmark problem in the simulation of turbulent flows. The main finding is that the power law relating the ratios at time $t$
to $T^*-t$ where $T^*$ is the peak enstrophy time is of a form that allows the machinery of dynamic interpolation-sparseness 
to produce a lower bound on the radius of spatial 
analyticity sufficient to overcome an upper bound on the scale of sparseness of the super-level
sets in view. As a consequence, the mechanism of turbulent dissipation engages via the harmonic measure maximum principle, 
furnishing a rigorous explanation for the subsequent slump of the enstrophy. 
This indicates that the 
higher-order derivative ratios -- which could be viewed as higher-order analogs of the classical Taylor and Kraichnan scales in turbulence phenomenology -- may
be reasonable identifiers of the peak of the energy dissipation rate.
\end{abstract}

\section{Introduction}

Higher-order derivative ratios emerged as \emph{key dynamic quantities} in the recent mathematical framework for the study of possible
singularity formation in the 3D Navier-Stokes and related systems based on a suitably defined scale of sparseness of the super-level 
sets of the higher-order spatial velocity derivatives \cite{Grujic2019, Grujic2020}. The underlying principle is that if the scale of sparseness 
falls below the scale of the
radius of spatial analyticity (which is a natural dissipation scale), the harmonic measure maximum principle will prevent any further
singularity buildup -- this is simply a rigorous mathematical description of the mechanism of turbulent dissipation \cite{Grujic2013}.
In the aforementioned framework, local-in-time dynamics of the chain of derivatives was studied and it was observed that monotonicity
of the portions of the chain, both `ascending' (higher order derivatives dominate lower order derivatives, consistent with turbulent buildup)
and `descending' (lower order derivatives dominate higher-order derivatives, consistent with turbulent decay), could be used to
improve upon the classical estimates. Tracking the ordered dynamics in time was accomplished precisely by studying local evolution
of the various derivative ratios.

\medskip

In particular, the utility of the ascending condition is in replacing the Gagliardo-Nirenberg interpolation inequalities in the higher-order Leibniz expansion of the 
nonlinearity. Since the derivative ratios are tracked dynamically, this was termed \emph{dynamic interpolation}. In the Navier-Stokes case, 
the general form of the ascending condition over a portion of the chain is as follows, 
\begin{equation}\label{w}
    \frac{\norm{D^ju(t)}{}^{\frac{1}{j+1}}}{\norm{D^ku(t)}{}^{\tfrac{1}{k+1}}}
     \lesssim (T-t)^{-\tfrac{1}{2} (\mu-1) \left(\tfrac{1}{j+ 1} - \tfrac{1}{k + 1}\right)} \dfrac{\left(j!\right)^{\frac{1}{j+1}}}{\left(k!\right)^{\tfrac{1}{k+1}}} , \qquad \ell \le j \le k
\end{equation}
for $t$ in some interval $(T-\epsilon, T)$ where $\mu$ is a suitable parameter and the norm is the $L^\infty$-norm. 
In the original work \cite{Grujic2019} the ascending building blocks
correspond to the case $\mu = 1$ (the ratios were bounded by constants), the time-weights allow for more flexibility.

\begin{remark}
Let us note that most of the scaling bounds on the derivative ratios to appear here (including the one above) include dimensional constants
that are hidden in "$\lesssim$" symbol.
\end{remark}

\medskip

In the case of the 3D hyper-dissipative Navier-Stokes system in the super-critical regime,
\begin{align}
&u_t+ (u\cdot\nabla) u = -(-\Delta)^\beta u-\nabla p               ,\\
&\textrm{div}\ u=0        
\end{align}
where the exponent $1<\beta <\frac{5}{4}$ measures the strength of the hyper-dissipation, the vector field $u$ is the velocity of the fluid and the scalar field
$p$ the pressure, the aforementioned methodology ruled out singularity formation in `turbulent regime' \cite{Grujic2020} which was 
determined by the condition
\[
    \frac{\norm{D^ku(t)}{}^{\frac{1}{k+1}}}{\norm{D^{2k}u(t)}{}^{\tfrac{1}{2k+1}}}
     \lesssim (T-t)^{-\frac{\beta-1}{2k+1}}
\]
(for sufficiently large $k$). In particular, this ruled out asymptotically self-similar blow-up, a prime candidate
for singularity formation in all super-critical hyper-dissipative Navier-Stokes models signaling that the infamous
Navier-Stokes super-criticality barrier may not be as resilient as it had been perceived.

\medskip

In the case of the Navier-Stokes system \emph{per se}, the methodology so far did not rule out any super-critical scenarios
(in particular, as pointed out in \cite{AlBr2022}, it did not rule out any super-critical generalized self-similar blow-ups) -- however -- it
did demonstrate asymptotic criticality (with respect to the order of the derivative) of the regularity problem
within the framework \cite{Grujic2019}.

\medskip

Given the above, it is natural to ask whether dynamical signatures of the higher-order derivative ratios are a reasonable identifier of the 
peak of the energy dissipation rate.
On one hand, we know from the rigorous theory that if the ratios scale (with respect to the local time-scale)
in a suitable way, the machinery of dynamic interpolation-sparseness will indeed trigger the mechanism of turbulent dissipation via the harmonic measure
maximum principle (the radius of spatial analyticity will dominate the scale of sparseness). On the other hand, so far there are no results 
(neither rigorous nor computational) indicating that the peak of the energy dissipation rate is preceded by
scaling of the ratios consistent with what is needed for dynamic interpolation-sparseness to successfully deploy.

\medskip

The main goal of the current work is to fill this void by examining scaling of the higher-order derivative ratios on a time interval leading to 
the enstrophy/energy dissipation rate peak in a computational simulation of the 3D Taylor-Green vortex (TGV) -- a benchmark problem in
simulating turbulent flows \cite{TaGr1937, BMONMF1983, DeB2013}. The following scaling transpired from our data analysis,

\[
    \frac{\norm{D^ku(t)}{}^{\frac{1}{k+1}}}{\norm{D^{2k}u(t)}{}^{\tfrac{1}{2k+1}}}
     \lesssim (T^*-t)^\frac{1}{k^\alpha}
\]
where $T^*$ is the peak enstrophy time (around 8.91 time units) and the exponent $\alpha \approx 0.89$. 
This was exciting to see since it will be argued in the theoretical part of the paper that any exponent less or equal to 1 would indeed
suffice to explain the subsequent slump of the enstrophy via dynamic interpolation-sparseness dissipation argument.

\medskip

The rest of the paper is organized as follows. Section 2 presents the computational results and the data analysis leading to the scaling law above,
Section 3 theoretical analysis prompted by the results in Section 2, and Section 4 conclusions.

\section{3D Taylor-Green vortex simulation and data analysis}

Consider the following form of the Navier-Stokes Equations (NSE)
\begin{align}
	u_t + (u \cdot \nabla) u - \nu \triangle u + \nabla p &= 0, &(x, t) \in \Omega \times (0, T), \label{momentum}\\
	\text{div} \, u &= 0, &(x, t) \in \Omega \times (0, T),\label{incomress} \\
	u(x, 0) &= u_0(x), &x \in \Omega\label{initial-con}
\end{align}\label{strong-form-problem}where $\Omega = [0, 1]^3$ is the unit cube with periodic boundary condition, and $T$ the simulation end time. The temporal and spatial variable are represented by $t$ and $x$, respectively. The viscosity is denoted by $\nu$, and $u_0$ is the initial condition. Denote by $\omega=\nabla\times u$ the vorticity and compute the following quantity of interest.

\smallskip

\noindent Enstrophy \begin{align}
\mathcal{E}(t)=\frac12\int_\Omega|\omega(x,t)|^2dx,
\end{align}
energy \begin{align}
E(t)=\frac12\int_{\Omega}|u (x,t)|^2 dx,
\end{align} and the derivative ratios
\begin{align}
R^k=\frac{\|D^{k}u(t)\|_{\infty}^{\frac{1}{k+1}}}{\|D^{2k}u(t)\|_{\infty}^{\frac{1}{2k+1}}}.
\end{align}   

Writing $u(x,t)=\sum\limits_{k}\hat{u}(k,t)e^{2\pi ik \cdot x}$, we convert \eqref{momentum}-\eqref{initial-con} to the Fourier space and arrive
at the following Initial Value Problem (IVP) 
\begin{align}
	\hat{u}_t  + \widehat{u\cdot\nabla u} +\nu|k|^2\hat{u}+ ik\hat{p} &= 0, \label{FT-1}\\
	k \cdot \hat{u} &= 0,\\
	\hat{u}(k,0)&=\hat{u}_0(k)\label{FT-2}
\end{align}
where $k$ is the wavevector. Applying Leray's orthogonal projector $$P_{k}=I-\frac{k\otimes k}{|k|^2}$$
to \eqref{FT-1}-\eqref{FT-2} yields
\begin{align}
	\hat{u}_t  +P_{k}( \widehat{u\cdot\nabla u}) +\nu|k|^2\hat{u}&= 0,\label{IVP-1}\\
	\hat{u}(k,0)&=\hat{u}_0(k).\label{IVP-2}
\end{align}

\begin{figure} [ht]
	\centering	
	\subfloat[]
	{\includegraphics[width=0.45\textwidth,height=0.36\textwidth]
		{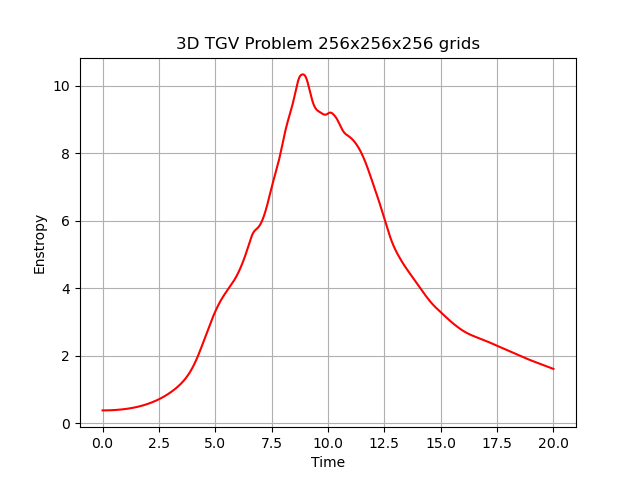}}\\\subfloat[]
	{\includegraphics[width=0.45\textwidth,height=0.36\textwidth]
		{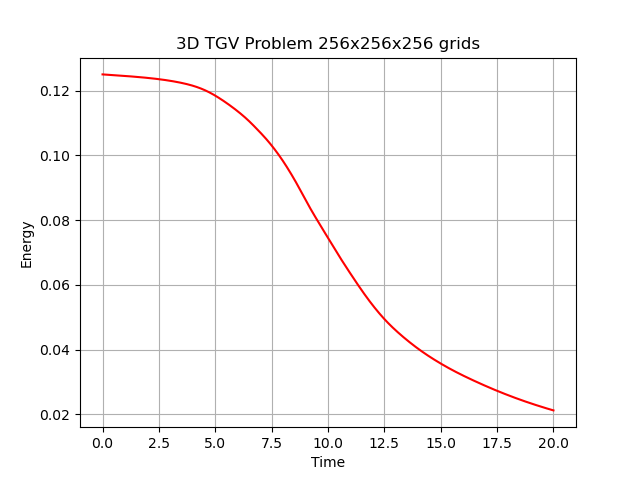}}
	\subfloat[]
	{\includegraphics[width=0.45\textwidth,height=0.36\textwidth]
		{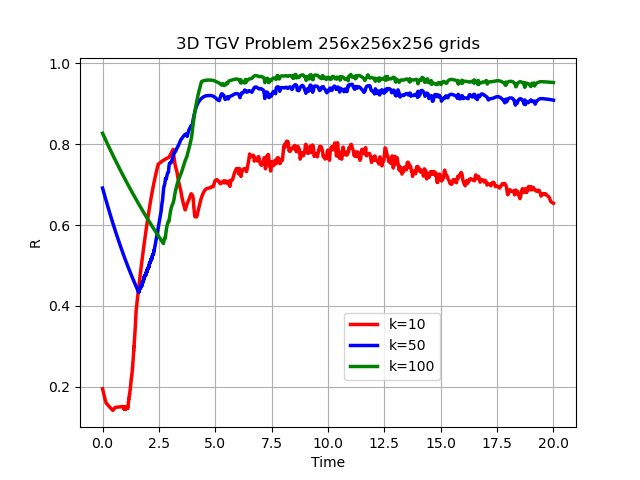}}\caption{\footnotesize{3D TGV problem with $256\times 256\times 256$ grid  (a) Enstrophy vs Time}, (b) Energy vs Time, and (c) Ratios vs Time graphs.}\label{3D-TG-256}
\end{figure}

Next we solve the benchmark 3D TGV problem \cite{TaGr1937, BMONMF1983, DeB2013}. To this end, we numerically solve the IVP in \eqref{IVP-1}–\eqref{IVP-2} using a fourth-order Runge–Kutta (RK4) time-stepping algorithm over the domain $\Omega=[0,1]^3$ and setting $\nu=1/1600$. Periodic boundary conditions are imposed, and the initial condition is given as follows \begin{align*}
	u_0(x) = \begin{pmatrix}
		\sin\left(2\pi x_1\right)\cos\left(2\pi x_2\right)\cos\left(2\pi x_3\right)\\\\
		-\cos\left(2\pi x_1\right)\sin\left(2\pi x_2\right)\cos\left(2\pi x_3\right)\\\\0
	\end{pmatrix}
\end{align*} 
where $x=(x_1,x_2,x_3)$. We discretize the domain into $256^3$ grid points, use a pseudo-spectral solver and run the simulation until $T=20$ with a uniform time-step size $\Delta t=0.001$.  We plot the Enstrophy vs. Time, and Energy vs. Time curves in Figure \ref{3D-TG-256}(a), and Figure \ref{3D-TG-256} (b), respectively; the plots are consistent with the classical benchmarks (cf. \cite{BMONMF1983, DeB2013}). In addition, we plot several Ratios vs. Time graphs in Figure \ref{3D-TG-256} (c) -- notice distinctive behavior of the ratios in the initial (laminar, vortex stretching intensifies gradients, $t \approx 0-4$), transitional (vortex sheets form, roll up and start breaking down, 
$t \approx 4-7/8$) and fully developed turbulence (starting at $t \approx 8$) regimes.

\medskip

From Figure \ref{3D-TG-256}(a), we estimate $T^*=8.91$ -- a critical value where the maximum enstrophy occurs -- and define $\beta=T^*-t$. The first step in our data analysis is to estimate the best $\gamma_k\in\mathbb{R}$ so that a power-type relationship $R^k = \beta^{\gamma_k}$ holds near $T^*$ and for a range of differential orders k. Using the available data $(\beta_i, R_i^k)$ where $\beta_i=\beta(t_i)$,
$R_i^k=R^k(t_i)$ and $t_i, i=1,2, \dots n$ are sourced from a unit-time interval preceding the enstrophy peak ($n \approx 10^3$), we find $\gamma_k$ in a least square sense by solving an overdetermined linear system via normal equation as below: 
 \begin{align}
	\gamma_k=\frac{\sum\limits_{i=1}^n\ln(\beta_i)\ln(R_i^k)}{\sum\limits_{i=1}^n\ln(\beta_i)^2}\label{gamma}
\end{align}
for each $k=k_1,k_2,\cdots,k_m$ (here, we choose $k_j=5 j, j = 1, 2,  \dots, 20$).

\medskip

\begin{figure} [ht]
	\centering	
	\subfloat[]{\includegraphics[width=0.45\textwidth,height=0.3\textwidth]
		{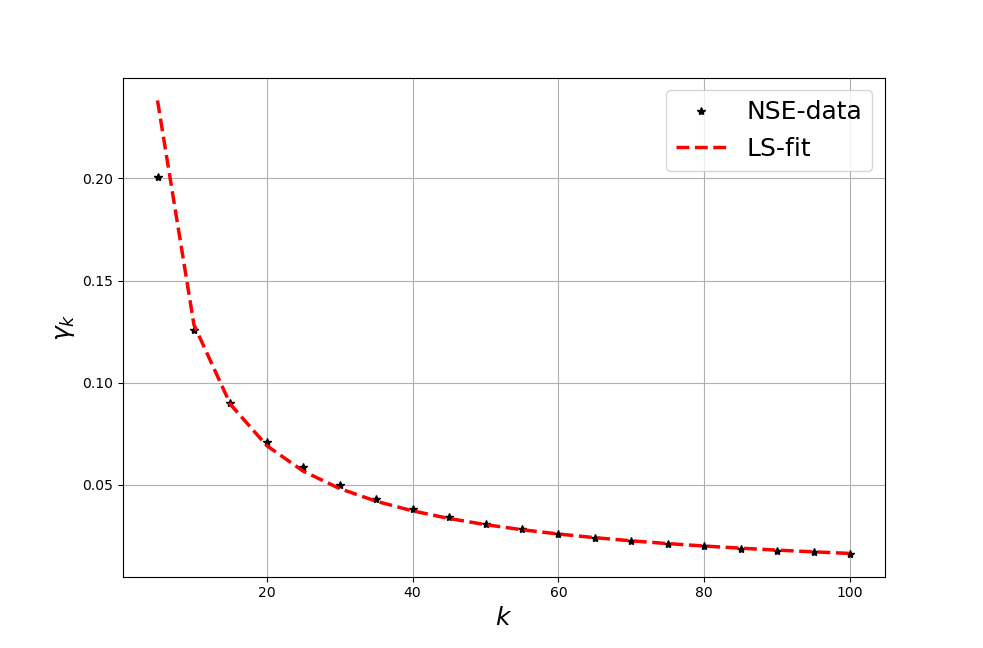}}
	\subfloat[]{\includegraphics[width=0.45\textwidth,height=0.3\textwidth]
		{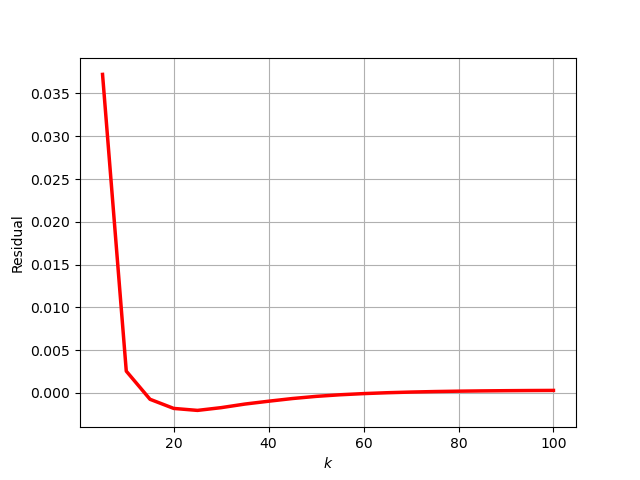}}
	\caption{\footnotesize{3D TGV problem with $256\times 256\times 256$ grid: (a) Least square estimate of $\gamma_k$ showing with  their fitted curve, and (b) Residuals vs. $k$.}}\label{gamma-fitted}
\end{figure}

Next we  plot the data $(\gamma_{k_1},\gamma_{k_2},\cdots,\gamma_{k_m})$ vs. $(k_1,k_2,\cdots,k_m)$ in Figure \ref{gamma-fitted}(a) with legend `NSE-data'. In the second step of our data analysis, we fit a model among the data sets $(k_1,k_2,\cdots,k_m)$ and $(\gamma_{k_1},\gamma_{k_2},\cdots,\gamma_{k_m})$ using Least-Square data fitting (LS-fit) model of the form 
\begin{align}
	\gamma_k=\frac{1}{k^a}
\end{align}
for some $a\in\mathbb{R}$. Following a similar approach as in \eqref{gamma}, we obtain the least-square estimate of $a=0.89$; that is, the best fit model is 
\begin{align}
	\gamma_k=\frac{1}{k^{0.89}}.
\end{align}

We conclude by plotting the LS-fit model outcomes in Figure \ref{gamma-fitted}(a) and the model residuals vs. $k$ in Figure \ref{gamma-fitted}(b).  

\clearpage

\section{Interpreting empirical scaling of the derivative ratios within dynamic interpolation-sparseness framework}

The goal of this section is to interpret the computational results presented in the previous section within the theoretical framework
of the dynamic interpolation-sparseness avenue to turbulent dissipation. Recall that the Navier-Stokes equations read

\begin{align}
&u_t+ (u\cdot\nabla) u = \Delta u-\nabla p               ,\\
&\textrm{div}\ u=0        
\end{align}
where the vector field $u$ is the velocity and scalar field $p$ the pressure (the viscosity is rescaled to 1 and the external force is
set to 0). Symmetrizing the nonlinearity via the incompressibility condition, the system can be written as follows,

\begin{align}\label{nse}
&u_t - \Delta u =   - \text{div} [u \otimes u] -\nabla p               ,\\
&\textrm{div}\ u=0.        
\end{align}

\medskip

Within this framework, turbulent dissipation will engage as soon as -- at some differential level -- the scale of the 
radius of spatial analyticity dominates the scale of sparseness of the super-level set of the maximal component of the derivative.
Recall that for a finite energy solution emanating from a smooth initial datum, an upper bound on 
the level-$k$ scale of sparseness is 
given by 
\[
r_k =  \|D^k u\|^{-\frac{1}{k+\frac{3}{2}}}
\]
\cite{Grujic2019}; hence if one could obtain a lower bound on the level-$k$ scale of the analyticity radius of
the form
\[
\rho_k =  \|D^k u\|^{-\frac{1}{(1+\epsilon_k)(k+1)}}
\]
where $\epsilon_k \ge \epsilon_*>0$ for $k$ large enough, $\rho_k$ would eventually surpass $r_k$ and dissipation would commence
(in the case of an ascending chain where the ratios are dominated by suitable constants, a lower bound on the 
scale of the analyticity radius emerged as 
\[
\|D^k u\|^{-\frac{1}{k+1}}
\]
\cite{Grujic2019}, meeting the scale of sparseness only asymptotically (in $k$).

\medskip

In what follows, we will sketch an argument revealing how the scaling law of the type that transpired in the previous section,
\[
    \frac{\norm{D^ku(t)}{}^{\frac{1}{k+1}}}{\norm{D^{2k}u(t)}{}^{\tfrac{1}{2k+1}}}
     \lesssim (T^*-t)^\frac{1}{k^\alpha}
\]
will indeed yield $\epsilon_k \ge \epsilon_* > 0$ as long as $\alpha \le 1$. Since our data analysis sourced from a fully
developed turbulence phase of the 3D Taylor-Green vortex flow leading to the enstrophy peak produced $\alpha \approx 0.89$,
this provides the first evidence that the higher-order derivative ratios may indeed be reasonable identifiers of the peak of the
energy dissipation rate.

\medskip

Denoting the heat semigroup by $S(t)$, Duhamel's principle yields the following form of (\ref{nse}),

\begin{equation}\label{duhamel}
 u(t)=S(t)u_0-\int_0^t S(t-s) (\text{div}[u\otimes u]+\nabla p)(s) \, ds.
\end{equation}

In this framework, level-$k$ dynamics is realized via evolution of spatial derivatives of order $k$. To this end,
taking $D^k$ of (\ref{duhamel}) yields

\begin{equation}\label{kduhamel}
 D^ku(t)=S(t)D^ku_0-\int_0^t \nabla \cdot S(t-s) (D^k[u\otimes u]+\nabla D^kp)(s) \, ds.
\end{equation}

A method for obtaining a lower bound on the radius of spatial analyticity for the Navier-Stokes equations in $L^q$-spaces 
pioneered in \cite{GrKu1998} is based on complexifying the equations, splitting the resulting system in real and imaginary parts and observing
that -- locally in time -- the radius of spatial analyticity grows at least like $\sqrt t$. This produces the maximal
lower bound on the analyticity radius comparable to the square root of the time of local existence of the 
solution. For a proper $q$ the pressure term is handled simply by the Calderon-Zygmund theorem since 
\[
 p = - (\triangle^{-1} \circ \partial_i\partial_j) u^iu^j
\]
and the operator  $\triangle^{-1} \circ \partial_i\partial_j$ is a composition of Riesz transforms. In the case $q=\infty$, the 
case we are interested here, handling the pressure is a bit more delicate -- one uses the fact  that the solution of the nonhomogeneous
heat equation with the right-hand side being a divergence of a $BMO$ function belongs to $L^\infty$ in conjunction with the 
Calderon-Zygmund theorem in $BMO$. In either case the bound on the pressure term is seamlessly incorporated in the bound on the nonlinear
term. In addition, it turns out that the local-in-time existence of the complexified system is comparable to the local-in-time existence of the
Navier-Stokes equations in the real space.

\medskip

Let us start by considering the nonlinear term. Using the Leibniz rule,
\[
 D^k[u \otimes u] = \sum_{j=0}^k \binom{k}{j} D^ju \otimes D^{k-j}u.
\]

Since the goal of this section is to pinpoint the role of the scaling signature of the derivative ratios obtained in the previous
section in improving the lower bound on the radius of spatial analyticity without fully engaging in the technical intricacies
of the rigorous argument,
we are going to consider the case of an even differential order $2k$ and retain only the symmetric part in the Leibniz
expansion. More precisely, we consider the following problem

\begin{equation}\label{2kduhamel}
 D^{2k}u(t)=S(t)D^{2k}u_0-\int_0^t \nabla \cdot S(t-s) ([D^ku \otimes D^ku]+\nabla D^{2k}p)(s) \, ds
\end{equation}
(in the actual proof, one would use a $\frac{1}{k^\alpha}$-version of the general weights displayed in (\ref{w}); this would cover the 
full range of indices appearing in the Leibniz expansion).

\medskip

In what follows, we are interested in local-in-time existence of the above system by performing a fixed point algorithm in
$L^\infty$. Since the estimate on the pressure term can be absorbed in the estimate on the nonlinear term (see the above
discussion), we are left to estimate
\[
 \Bigl\| \int_0^t \nabla \cdot S(t-s) (D^ku \otimes D^ku)(s) \, ds \Bigr\|.
\]

Assuming the scaling signature of the derivative ratios of the form
\[
    \frac{\norm{D^ku(s)}{}^{\frac{1}{k+1}}}{\norm{D^{2k}u(s)}{}^{\tfrac{1}{2k+1}}}
     \lesssim (t-s)^\frac{1}{k^\alpha}
\]
one obtains
\[
 \|D^ku(s)\|^2 \lesssim (t-s)^\frac{2(k+1)}{k^\alpha} \|D^{2k}u(s)\|^\frac{2(k+1)}{2k+1}
\]
which -- in turn -- yields (in conjunction with the estimate on the heat kernel)
\[
  \Bigl\| \int_0^t \nabla \cdot S(t-s) (D^ku \otimes D^ku)(s) \, ds \Bigr\| \lesssim t^{\frac{1}{2}+\frac{2(k+1}{k^\alpha}}  \sup_{s \in (0, t)} 
  \|D^{2k}u(s)\|^\frac{2(k+1)}{2k+1}.
\]
If we denote the maximal existence time guaranteed by the above estimates by $T$, then
\[
 T^{\frac{1}{2}+\frac{2(k+1)}{k^\alpha}} \|D^{2k}u_0\|^\frac{1}{2k+1} \lesssim 1
\]
which leads to the following lower bound on the scale of the analyticity radius,
\[
\rho_{2k} =  \|D^{2k} u\|^{-\frac{1}{(1+\epsilon_{2k})(2k+1)}}
\]
where
\[
 \epsilon_{2k}=\frac{4(k+1)}{k^\alpha}.
\]
Since $\epsilon_{2k}$ stays bounded away from zero for any $\alpha \le 1$, it transpires that $\rho_{2k}$ will dominate
the upper bound on the scale of sparseness
\[
r_{2k} =  \|D^{2k} u\|^{-\frac{1}{2k+\frac{3}{2}}}
\]
in this scaling regime, making it possible for dissipation to engage via the harmonic measure maximum principle.

\begin{remark}
In order to make the above sketch rigorous one also needs to propagate the scaling assumptions/properties of the derivative 
ratios locally-in-time and into the complex space; this introduces an additional level of technical complexity \cite{Grujic2019, Grujic2020}
and is outside of scope of this note. 
\end{remark}

\section{Conclusion}

The goal of this note was to examine the scaling properties of the derivative ratios on a time interval
leading to the peak of the enstrophy/energy dissipation rate in the case of the 3D Taylor-Green vortex -- a benchmark problem in simulations of turbulent
flows. Recent work \cite{Grujic2019, Grujic2020} showed that certain monotonicity properties of the portions of the chain of derivatives,
equivalently certain scaling properties of the various derivative ratios, could successfully engage the mechanism of turbulent dissipation
via the harmonic measure maximum principle within the dynamic interpolation-sparseness mathematical framework. 
However, there were no results on whether the scaling signatures of these
ratios might also be an identifier of the peak of the energy dissipation rate, rather than simply a predictor.
This note is the first contribution in this direction, revealing that in the case of the 3D Taylor-Green vortex the scaling signatures of
the ratios sourced from the computational data in the regime of interest are indeed capable of explaining the slump past the 
enstrophy/energy dissipation peak.

\def\cprime{$'$}

\end{document}